\documentclass[oneside]{amsart}

\usepackage{control/packages}
\usepackage{control/docaesthetic}
\usepackage{control/environments}
\usepackage{control/macros}
\usepackage{control/paperdata}

\usepackage{subfiles}

\usepackage[style=numeric]{biblatex}
\begin{document}

\maketitle

%

\begin{abstract}
We introduce a generalization of the product expansion of a finite semigroup.
As an application, we provide an alternative proof of the decidability of pointlike sets for pseudovarieties consisting of semigroups whose subgroups all belong to a given decidable pseudovariety of groups.
\end{abstract}

%

\tableofcontents

%

\section{Introduction}

The purpose of this paper is twofold.
First, we develop a generalization of the product expansion of a finite semigroup. 
This construction was introduced by the first author in \cite{PRODEXP}, where it was used to obtain a simplified proof that aperiodic pointlikes are decidable.\footnotemark
\footnotetext{The original proof of this result is due to the first author in \cite{PL-APERIODIC-KARSTEN}.} 
As such, it should not come as a surprise that our generalized construction affords a generalized result---to be precise, we obtain an alternative proof of Steinberg and van Gool’s result from \cite{PL-VARDETGRP-BEN}, which is that the pseudovariety
    \[
        \pvGH = \setst{S \in \FSGP}{\text{all subgroups of $S$ belong to $\pvH$}}
    \]
has decidable pointlikes whenever $\pvH$ is a pseudovariety of groups in which membership is decidable. 

The second purpose of this paper is to serve as an illustration of the authors' general framework for pointlike sets as developed in \cite{GENTHEORYPL}. 
That being said, this paper is intended to be as self-contained as possible---all notions and results from \cite{GENTHEORYPL} which are needed in this paper are reiterated here (albeit without proof in the case of results). 

\subsection{Organization of paper}
Section \ref{section:prelim} establishes notation and reviews various standard concepts and results which are needed later.
Section \ref{section:framework} covers a number of concepts and results from \cite{GENTHEORYPL} which provide the framework for our approach.
This material is then used in Section \ref{section:lowerbound} to define our ``candidate'' for $\bbG \pvH$-pointlikes as well as to show that this candidate provides a lower bound.
Section \ref{section:flows} considers flows on finite automata as a method for constructive proofs of upper bounds for pointlikes. 
The material covered in Section \ref{section:blowup}---namely, blowup operators---is the final piece required to develop the titular product expansion in Section \ref{section:prodexp}.\footnotemark
\footnotetext{It should be noted that the the product expansion is constructed with respect to an arbitrarily chosen collection of ``setup data'', and hence the construction is non-functorial.}
In Section \ref{section:cascades}, we review notions relating to wreath products of transformation semigroups---in particular, we review the so-called ``Zeiger property'' and its relation to subgroups.
This leads into Section \ref{section:properties}, wherein we show that subgroups of the transition semigroup of the product expansion are bound in an advantageous manner by the setup data. 
Finally, Section \ref{section:upperbound} uses the product expansion to establish the upper bound, thereby proving the main result.

%

%

\section{Preliminaries}     \label{section:prelim}

\subsection{}
Familiarity with finite semigroup theory is assumed; for reference, the reader is directed to \cite{QTHEORY}.

\begin{notation}
Let $S$ be a finite semigroup.
\begin{itemize}
    \item Write $S^I$ for the semigroup obtained by adjoining a new element $I$ to $S$ and defining $xI = Ix = x$ for all $x \in S$.
    \item Write $S^0$ for the semigroup obtained by adjoining a new element $0$ to $S$ and defining $x0 = 0x = 0$ for all $x \in S$.
    \item Given $x \in S$, let $x^\omega$ denote the unique idempotent generated by $x$.
    \item Green's equivalence relations are denoted by $\GR$, $\GL$, $\GH$, and $\GJ$; 
    moreover, the various Green's equivalence classes of $x \in S$ are denoted by $\RCL{x}$, $\LCL{x}$, $\HCL{x}$, and $\JCL{x}$, respectively.
\end{itemize}
\end{notation}

\subsection{Activators}
If $J$ is a $\GJ$-class of a finite semigroup $S$, then the set
    \[
    A_L(J) \;=\; \setst{a \in S^I}{aJ \cap J \neq \varnothing}
    \]
is a union of $\GJ$-classes of $S^I$.
Moreover, $A_L(J)$ contains a unique $\lgJ$-minimal $\GJ$-class which is called the \textbf{left activator} of $J$ and which will be denoted by $\LACT(J)$.
Note that $J$ is regular if and only if $\LACT(J) = J$.

\begin{lem} \label{activatorlemma}
Let $J$ be a $\GJ$-class of $S$.
\begin{enumerate}
    \item For each $j \in J$, there exists $a \in \LACT(J)$ such that $a j = j$.
    In particular, there exists an idempotent $e \in \LACT(J)$ such that $e j = j$.
    \item Given $j \in J$ and $a \in \LACT(J)$ as in (1), 
        \[
            xj \slgJ j 
                \quad \Longleftrightarrow \quad
            xa \slgJ a
        \]
    for all elements $x \in S$.
\end{enumerate}
\end{lem}

\begin{proof}
See \cite[Lemma~2.9]{CPL-LOWERBOUNDS}.
\end{proof}

\begin{notation}
Given $x \in S$, write $\LACT(x) = \LACT( \JCL{x} )$. 
\end{notation}

\subsection{Relational morphisms}
A \textbf{relational morphism} $\rho : S \rlm T$ is an equivalence class of spans of the form
    \begin{equation*}
        \begin{tikzcd}
            \cdot 
            \arrow[r] 
            \arrow[d, two heads] 
        & 
            T 
        \\
            S 
        &
        \end{tikzcd}
    \end{equation*}
where the map to $S$ is a regular epimorphism,\footnotemark and where two such spans are equivalent if the natural maps from each apex to $S \times T$ have the same image.
\footnotetext{Here and throughout, ``regular epimorphism'' means ``surjective homomorphism''.}
This shared image is called the \textbf{graph} of $\rho$ and is denoted by $\Gamma(S , \rho, T)$.

\subsection{Pseudovarieties}
A \textbf{pseudovariety} is a class of finite semigroups which is closed under taking subsemigroups, quotients, and finite products of its members.

\subsection{Power semigroups}   \label{prelim:power}
Given a finite semigroup $S$, let $\Po(S)$ denote the semigroup of non-empty subsets of $S$ under the inherited operation given by
    \[
        X \cdot Y = \setst{xy}{x \in X, y \in Y}
    \]
for all non-empty subsets $X$ and $Y$ of $S$; also, let $\sing(S)$ denote the subsemigroup of $\Po(S)$ consisting of the singletons.

A morphism $\varphi : S \rightarrow T$ extends to a morphism
    \[
        \ext{\varphi} : \Po(S) \arw \Po(T)
            \quad \text{given by} \quad
        (X)\ext{\varphi} = \setst{(x)\varphi}{x \in X}.
    \]
Equipping the object map $\Po$ with this action on morphisms yields a functor
    \[
        \Po : \FSGP \arw \FSGP
    \]
which creates monomorphisms, regular epimorphisms, and isomorphisms.

\subsection{Pointlikes}
Given a finite semigroup $S$ and a pseudovariety $\pvV$, a non-empty subset $X \subseteq S$ is said to be \textbf{$\pvV$-pointlike} if for any relational morphism of the form $\rho : S \rlm V$ with $V \in \pvV$ there exists some element $v \in V$ for which $X \subseteq (v)\rho^\inv$.

The set of $\pvV$-pointlike subsets of $S$ is denoted by $\PV(S)$, and is in fact a subsemigroup of $\Po(S)$ which contains $\sing(S)$ and which is closed under taking non-empty subsets of its members.
Equipping this object map with the action on morphisms sending $\varphi : S \rightarrow T$ to the evident restriction of the extension described in \ref{prelim:power} yields a subfunctor
    \[
        \PV : \FSGP \arw \FSGP
    \]
of $\Po$ with the property that a finite semigroup $S$ belongs to $\pvV$ if and only if $\PV(S) = \sing(S)$.\footnotemark
\footnotetext{The notation for $\Po$ is due to it being the pointlikes functor for the trivial pseudovariety $\pvtriv$.}
Pointlike functors also create monomorphisms, regular epimorphisms, and isomorphisms.

%

%

\section{Review of framework}       \label{section:framework}

\subsection{}
The concepts considered in this section are covered in minimal detail; we discuss only the notions and results which are necessary for the rest of the paper.
For more detail, see \cite{GENTHEORYPL}.

\subsection{Semigroup complexes}
Let $S$ be a finite semigroup.
An \textbf{$S$-complex} $\cK$ is a subsemigroup of $\Po(S)$ which 
\begin{enumerate}
    \item contains $\sing(S)$ as a subsemigroup, and which
    \item is closed under taking non-empty subsets of its members, meaning that if $X \in \cK$ then any $Y \in \Po(S)$ for which $Y \subseteq X$ is also a member of $\cK$.
\end{enumerate}
The set of $S$-complexes, denoted by $\Com{S}$, is a complete lattice wherein the order is inclusion, the top and bottom are $\Po(S)$ and $\sing(S)$ respectively, the meet is intersection, and the join is given by computing the join in the lattice of subsemigroups of $\Po(S)$ then closing under taking non-empty subsets.

The \textit{raison d'{\^e}tre} for $S$-complexes is that $\PV(S)$ is an $S$-complex for any pseudovariety $\pvV$, and hence $\Com{S}$ serves as a ``local possibility space'' for pointlikes.

\begin{defn}
A \textbf{modulus} $\Lambda$ is a rule, generally written as 
    \[
        \Lambda = \modrl{S}{\Lambda_S},
    \]
which assigns to each finite semigroup $S$ a (possibly empty) set $\Lambda_S \subseteq \Po(S)$, and which satisfies the following axioms.
\begin{enumerate}
    \item If $\varphi : S \rightarrow T$ is a morphism, then for any $X \in \Lambda_S$ there exists some $\widetilde{X} \in \Lambda_T$ such that $(X)\ext{\varphi} \subseteq \widetilde{X}$.
    \item If $\varphi : S \surj T$ is a regular epimorphism, then for any $Y \in \Lambda_T$ there exists some $\widetilde{Y} \in \Lambda_S$ such that $(\widetilde{Y})\ext{\varphi} = Y$.
\end{enumerate}
\end{defn}

\begin{example}
Commonplace moduli include
\begin{enumerate}
    \item the \textit{subgroup} modulus
        \[
            \mathsf{Grp} = \modrl{S}{ \setst{G}{ \textnormal{$G$ is a subgroup of $S$} } };
        \]
    \item the three \textit{Green's} moduli
        \[
        \mathsf{RCl} = \modrl{S}{S / \GR},
            \quad
        \mathsf{LCl} = \modrl{S}{S / \GL},
            \quad \text{and} \quad
        \mathsf{JCl} = \modrl{S}{S / \GJ};
        \]
    \item the three ``principal'' moduli, including the \textit{principal right ideals} modulus
        \[
            \mathsf{PrinR}
                =
            \modrl{S}{ \setst{ x \cdot S^I }{x \in S} },
        \]
    the \textit{principal left ideals} modulus
        \[
            \mathsf{PrinL}
                =
            \modrl{S}{ \setst{S^I \cdot x}{x \in S} },
        \]
    and the \textit{principal (two-sided) ideals} modulus
        \[
            \mathsf{PrinJ}
                =
            \modrl{S}{ \setst{S^I \cdot x \cdot S^I}{x \in S} };
        \]
    \item and the \textit{set of idempotents} modulus 
        \[
            \mathsf{E} = \modrl{S}{\{E(S)\}}.
        \]
\end{enumerate}
\end{example}

\subsection{Constructing lower bounds for pointlikes}   \label{framework:lowerbound}
Moduli provide a uniform way to describe the construction of candidates for pointlikes.
Given a modulus $\Lambda$, the \textbf{$\Lambda$-construct} of $S \in \FSGP$ is given by
    \[
        \scC_\Lambda(S)
            =
        \bigcap \setst{ \cK \in \Com{S} }{ \textnormal{if $\cX \in \Lambda_\cK$, then $\bigcup \cX \in \cK$} };
    \]
that is, $\scC_\Lambda(S)$ is the minimal $S$-complex closed under unioning subsets assigned to it by the modulus $\Lambda$.
Equipping the object map $\scC_\Lambda$ with the action on morphisms sending $\varphi : S \rightarrow T$ to the extension
    \[
        \ext{\varphi} : \scC_\Lambda(S) \arw \scC_\Lambda(T)
            \quad \text{given by} \quad
        (X)\ext{\varphi} = \setst{(x)\varphi}{x \in X}
    \]
yields a functor which, moreover, admits a monad structure $(\scC_\Lambda , \sigma_\Lambda , \mu_\Lambda)$,
where the components of the unit $\sigma_\Lambda : \id{\FSGP} \Rightarrow \scC_\Lambda$ are the singleton embeddings
    \[
        \sigma_{\Lambda , S} 
            =
        \{-\} : S \longinj \scC_\Lambda(S)
            \quad \text{given by} \quad
        x \longmapsto \{x\}
    \]
and the components of the multiplication $\mu_\Lambda : \scC_\Lambda^2 \Rightarrow \scC_\Lambda$ are the union maps
    \[
        \mu_{\Lambda , S} 
            =
        \bigcup (-) : \scC_\Lambda^2(S) \longsurj \scC_\Lambda(S)
            \quad \text{given by} \quad
        \cX \longmapsto \bigcup_{X \in \cX} X
    \]
for every finite semigroup $S$.

The set of \textbf{points} of a modulus $\Lambda$ is given by
    \[
        \points{\Lambda}    
            \;=\;
        \setst{S \in \FSGP}{\Lambda_S \subseteq \sing(S)}.
    \]
It is straightforward to see that $\points{\Lambda}$ is a pseudovariety (see \cite[Proposition~9.7]{GENTHEORYPL}).
Moreover, $S \in \points{\Lambda}$ if and only if $\scC_\Lambda(S) = \sing(S)$ for any $S \in \FSGP$.

The utility of this fact comes from \cite[Theorem~9.12]{GENTHEORYPL}, which states that if $\Lambda$ is a modulus with $\points{\Lambda} = \pvV$, then $\scC_\Lambda(S) \subseteq \PV(S)$ for all $S \in \FSGP$.
Consequently, constructing lower bounds for pointlikes is reducible to constructing moduli with the relevant pseudovariety of points.

\subsection{Nerves of relational morphisms}
The basic justification for this lower bound is as follows.
The \textbf{nerve} of a relational morphism $\rho : S \rlm T$ is the $S$-complex
    \[
        \Nrv(S , \rho, T)
            =
        \setst{ X \in \Po(S)}{\text{$X \subseteq (t)\rho^\inv$ for some $t \in T$}};
    \]
and in this situation we say that $\rho$ \textbf{computes} $\Nrv(S, \rho, T)$.

If $\Lambda$ is a modulus with $\points{\Lambda} = \pvV$ and $\rho : S \rlm V$ is a relational morphism with $V \in \pvV$, then it is easily seen that any $X \in \Lambda_S$ must be a member of $\Nrv(S, \rho, V)$.
From here, the well-known fact that there always exists a relational morphism of this form which computes $\PV(S)$ yields the desired bound.

%

%

\section{Candidate and lower bound}      \label{section:lowerbound}

\subsection{}
For the remainder of the paper, let $\pvH$ be a pseudovariety of groups.
Recall from the introduction that 
    \[
        \pvGH = \setst{S \in \FSGP}{\text{all subgroups of $S$ belong to $\pvH$}}.
    \]

\subsection{Candidate}      \label{lowerbound:lb}
Define a modulus $\Lambda_\pvGH$ by setting
    \[
        \Lambda_\pvGH
            =
        \modrl{S}{ \bigcup \setst{ \PL{\pvH}(G)}{\text{ $G$ is a subgroup of $S$ } } },
    \]
and denote the $\Lambda_\pvGH$-construct by $\CGH$; that is, $\CGH(S)$ is the minimal $S$-complex closed under unioning $\pvH$-pointlike subsets of its subgroups for any $S \in \FSGP$.

It is clear that $\points{\Lambda_\pvGH} = \pvGH$, and it follows that $\CGH(S) \subseteq \PGH(S)$ for any finite semigroup $S$ by \ref{framework:lowerbound}.

\subsection{}   \label{lowerbound:Hkernel}
If $\rho : G \rlm H$ is a relational morphism of groups, then $(1_H)\rho^\inv$ is a normal subgroup of $G$ for which $G / (1_H)\rho^\inv \cong \Img(\rho)$ (see \cite[Lemma~2.14]{PIN-RELMORLANG}), and so
    \[
        \Nrv(G, \rho, H) 
            =
        \setst{ X \in \Po(G) }{\text{$X \subseteq g \cdot (1_H)\rho^\inv$ for some $g \in G$}}.
    \]
Consequently, $\PL{\pvH}(G)$ is computable by a morphism for any finite group $G$.
In fact, there is a canonical choice of such a morphism---since the set of normal subgroups $N$ of $G$ for which the quotient $G/N$ belongs to $\pvH$ is closed under intersection, there is a unique minimal normal subgroup of $G$ with this property; this subgroup is called the \textbf{$\pvH$-kernel} of $G$ and is denoted by $K_{\pvH, G}$.
The quotient map $G \surj G / K_{\pvH , G}$---through which any homomorphism $G \rightarrow H$ with $H \in \pvH$ factors uniquely---is easily seen to compute $\PL{\pvH}(G)$ as its nerve.

\subsection{}
The observations in \ref{lowerbound:Hkernel} imply that
    \[
        \Lambda_\pvGH
            =
        \modrl{S}{\setst{K_{\pvH , G}}{\text{$G$ is a subgroup of $S$}}}.
    \]
This alternative characterization of $\Lambda_\pvGH$---which reflects the approach of \cite{PL-VARDETGRP-BEN}---will be highly relevant again later on in this paper.


%

\section{Flows and constructive upper bounds}      \label{section:flows}

\begin{defn}
An \textbf{automaton} is the data of a tuple $\cA = (\Sigma, Q, \INIT, \INP)$, where
\begin{enumerate}
    \item $\Sigma$ is a finite set of \textbf{input symbols},
    \item $Q$ is a finite set of \textbf{states},
    \item $\INIT \in Q$ is a distinguished \textbf{initial state}, and
    \item $(- \INP -) : Q \times \Sigma \rightarrow Q$ is a set function called the \textbf{transition function}.
\end{enumerate}
\end{defn}

\begin{notation}
\hfill
\begin{enumerate}
    \item The set of (non-empty) strings over $\Sigma$ is denoted by $\Sigma^+$.
    \item Write the cumulative action of a string $\sigma = a_1 a_2 \cdots a_n \in \Sigma^+$ on a state $q \in Q$ as $q \INP \sigma$; i.e,
        \[
            (- \INP \sigma) 
                \quad = \quad
            ( \dots ( (- \INP a_1) \INP a_2) \dots ) \INP a_n.
        \]
\end{enumerate}
\end{notation}

\begin{defn}
The \textbf{transition semigroup} of an automaton $\cA$ is the semigroup $\scT_\cA$ which is generated by the functions $(- \INP a) :  Q \rightarrow Q$ induced by each $a \in \Sigma$.
\end{defn}

\begin{example}
The \textbf{Cayley automaton} of a finite semigroup $S$ is the automaton
    \[
        \Cay(S) = (S, S^I, I, \cdot).
    \]
That is, the states of $\Cay(S)$ are the elements of $S$ along with an adjoined identity which behaves as the initial state, and the transition function is given by right multiplication.
It is easily verified that the transition semigroup of $\Cay(S)$ is isomorphic to $S$.
\end{example}

\begin{defn}
A \textbf{flow automaton} is a triple $(S, \cA, \flow{-}{})$ consisting of a finite semigroup $S$, an automaton $\cA = (S, Q, \INIT, \INP)$, and a set function
    \[
        \flow{-}{} : Q \setminus \{ \INIT \} \arw \Po(S),
    \]
which is the nominal \textit{flow} and which satisfies
    \[
    s \in \flow{\INIT \INP s}{}
    \qquad \textnormal{and} \qquad
    \flow{q}{} \cdot \{s\} \;\subseteq\; \flow{q \INP s}{}
    \]
for all $s \in S$ and every non-initial state $q \in Q \setminus \{\INIT\}$.
\end{defn}

\begin{lem} \label{automata:flow:lem:words}
If $(S, \cA, \flow{-}{})$ is a flow automaton, then $[\sigma]_S \in \flow{\INIT \INP \sigma}{}$ for all $\sigma \in S^+$.
\end{lem}

\begin{proof}
This is easily seen via induction on the length of $\sigma$.
\end{proof}

\begin{defn}
The \textbf{cover complex} of a flow automaton $(S, \cA, \flow{-}{})$ is the $S$-complex given by
    \[
        \Cov(S , \cA, \flow{-}{}) 
            = 
        \setst{X \in \Po(S)}{\text{$X \subseteq \flow{q}{}$ for some $q \in Q \setminus \{\INIT\}$}}
    \]
\end{defn}

\begin{prop} \label{flowequivalence}
If $S$ is a finite semigroup., then
    \[
        \PV(S)
        \;=\;
        \bigcap \setst{\Cov(S , \cA, \flow{-}{})}{\scT_\cA \in \pvV}
    \]
for any pseudovariety $\pvV$.
\end{prop}

\begin{proof}
See \cite[Proposition~2.5]{PL-VARDETGRP-BEN}.
\end{proof}

%

%

\section{Blowup operators}      \label{section:blowup}

\subsection{}
Fix a finite semigroup $S$ for the remainder of this section as well as Section \ref{section:prodexp}.

\begin{defn}
A \textbf{preblowup operator} on an $S$-complex $\cK$ is a set function
    \[
        \beta : \cK \arw \cK
    \]
which satisfies the following axioms.
\begin{enumerate}
    \item $X \subseteq (X)\beta$ for all $X \in \cK$;
    \item $X \gL Y$ implies $(X)\beta \gL (Y)\beta$ for any $X , Y \in \cK$; and
    \item either $(X)\beta = X$ or $(X)\beta \slgH X$ for each $X \in \cK$.
\end{enumerate}
An idempotent preblowup operator is called a \textbf{blowup operator}.
\end{defn}

\subsection{}
The set of preblowup operators on $\cK \in \Com{S}$ is closed under composition and is therefore a finite semigroup; consequently, any preblowup operator $\beta$ induces a blowup operator $\beta^\omega$.

\subsection{Multipliers for preblowup operators}
A \textbf{right multiplier} for a preblowup operator $\beta$ on $\cK$ is a set function
    \[
        r_\ast : \cK \arw \cK^I
            \quad \text{written} \quad
        X \longmapsto r_X
    \]
such that $X \cdot r_X = (X)\beta$ for all $X \in \cK$.

Dually, a \textbf{left multiplier} for $\beta$ is a set function
    \[
        \ell_\ast : \cK \arw \cK^I
            \quad \text{written} \quad
        X \longmapsto \ell_X
    \]
such that $\ell_X \cdot X = (X)\beta$ for all $X \in \cK$.

\subsection{Properties of preblowup operators}
Let $\beta$ be a preblowup operator on $\cK$, and let $r_\ast$ and $\ell_\ast$ be right and left multipliers for $\beta$, respectively.
The following properties are easily verified.
\begin{enumerate}
    \item If $X \gJ (X)\beta$, then $X = (X)\beta$.
    \item If $X \gJ Y$, then $(X)\beta \gJ (Y)\beta$.
    \item For both of the claims above, the analogous statement is true for $\GR$, $\GL$, and $\GH$ as well.
    \item Let $A, B, C, D \in \cK$. Then
        \begin{align*}
            AX = BX 
                \quad & \Longrightarrow \quad 
            A \cdot (X)\beta \;=\; B \cdot (X)\beta,
        \\
            XC = XD  
                \quad & \Longrightarrow \quad 
            (X)\beta \cdot C \;=\; (X)\beta \cdot D.
        \end{align*}
    More generally,
        \begin{align*}
            AX \gL BX 
                \quad & \Longrightarrow \quad 
            A \cdot (X)\beta \;\gL\; B \cdot (X)\beta,
        \\
            XC \gR XD  
                \quad & \Longrightarrow \quad 
            (X)\beta \cdot C \;\gR\; (X)\beta \cdot D.
        \end{align*}
    \item If $X \lgL Y$, then $X \subseteq X \cdot r_Y$. 
    Moreover, if $(Y)\beta = Y$, then $X =  X \cdot r_Y$.
    \item If $X \lgR Y$, then $X \subseteq \ell_Y \cdot X$. 
    Moreover, if $(Y)\beta = Y$, then $X =  \ell_Y \cdot X$.
\end{enumerate}
It follows from (1) that if $J$ is a $\GJ$-class of $\cK$ then either $(X)\beta = X$ for all $X \in J$ or every element of $J$ falls in both $\GL$ and $\GR$ to some lower $\GJ$ class.
The set of $\GJ$-classes which are fixed by $\beta$ is called the \textbf{$\beta$-core} of $\cK$ and is denoted by $\Core(\cK ; \beta)$, i.e.,
    \[
        \Core(\cK ; \beta)
            =
        \setst{J \in \cK / \GJ}{\text{$(X)\beta = X$ for all $X \in J$}}.
    \]
Note that the minimal $\GJ$-class of $\cK$ always belongs to the $\beta$-core of $\cK$.

%

%

\section{Construction of the product expansion}      \label{section:prodexp}

\subsection{}
The product expansion is a construction which, roughly speaking, ``blows up unwanted groups'' via a blowup operator (hence the name) by disassembling the $\GJ$-structure of an $S$-complex, throwing out $\GJ$-classes which aren't fixed by said blowup operator, then ``integrating'' the local parts into a coherent whole.

\subsection{Setup data}     \label{prodexp:setupdata}
A product expansion on $S$ is defined with respect to a tuple 
    \[
        \mathfrak{d} 
            =
        \left( \cK , \beta , E_\ast , J_\ast \right),
    \]
referred to as \textbf{setup data}, where
\begin{enumerate}
    \item $\cK$ is an $S$-complex;
    \item $\beta$ is a blowup operator on $\cK$;
    \item $E_\ast$ is a set function
        \[
            E_\ast : \cK \arw \cK^I
                \quad \text{written} \quad
            X \longmapsto E_X
        \]
    assigning to each $X \in \cK$ an idempotent $E_X \in \LACT(X)$---where $\LACT(X)$ is the left activator of the $\GJ$-class of $X$---such that $E_X X = X$; and\footnotemark
    \item $J_\ast$ is an antitone bijection 
        \[
            J_\ast : (\{1 , 2, \dots, N\} , \leq) \arw (\Core(\cK ; \beta) , \lgJ)
        \]
    where $N$ is the size of $\Core(\cK ; \beta)$;
    that is, it is an indexing of the $\beta$-core $\GJ$-classes of $\cK$ with the property that
        \[
             J_k \slgJ J_\ell 
                \quad \Longrightarrow \quad
            k > \ell
        \]
    for all $J_k , J_\ell \in \Core(\cK ; \beta)$.
\end{enumerate}
\footnotetext{The existence of such idempotents is guaranteed by Lemma \ref{activatorlemma}.}

\begin{notation}        \label{prodexp:notationfixed}
For the remainder of the construction we fix a tuple $\mathfrak{d}$ of setup data as in \ref{prodexp:setupdata}.
Moreover, 
\begin{itemize}
    \item write $\clos{X} = (X)\beta = (X)\beta^\omega$ for each $X \in \cK$;
    \item for each $J_k \in \Core(\cK ; \beta)$, let $d_k$ denote the length of the longest strict $\GJ$-chain---\textit{in $\cK$!}---descending from some element of $J_k$;
    \item $\cB(J_k) = \bigcup \setst{J \in \cK / \GJ}{J \slgJ J_k}$;
    \item $\Below(J_k) = \bigcup \setst{J_i \in \Core(\cK ; \beta)}{k > i}$; and
    \item $\Inputs(J_k) = \Below(J_k) \cup J_k \cup \{ I \}$.
\end{itemize}
\end{notation}

\subsection{}
The product expansion of $S$ with respect to the setup data $\mathfrak{d}$ is an automaton which is constructed as follows.
For each $\beta$-core $\GJ$-class $J_k$ we define
\begin{itemize}
    \item an automaton $(\Inputs(J_k), Q(J_k), \bullet, \odot)$,
    \item an ``output'' function $\prodout{-}{-} : Q(J_k) \times \Inputs(J_k) \rightarrow \cB(J_k) \cup \{ I \}$, and
    \item a ``value'' map $\prodval{-} : Q(J_k) \rightarrow \cK^I$,
\end{itemize}
each of whose precise details depend on whether $J_k$ is regular or null.
These automata are then ``composed'' in series via the order provided by $J_\ast$, where the input to each ``lower'' coordinate is obtained by applying $\beta$ to the output of the action above it.

\subsection{}
First, although the regular and null constructions are in most aspects quite different, the following definitions coincide for any $\beta$-core $\GJ$-class $J_k$:
\begin{itemize}
    \item $\prodval{\bullet} = I$;
    \item if $q \in Q(J_k)$, then
        \[
            q \odot I = q 
                \quad \text{and} \quad
            \prodout{q}{I} = I;
        \]
    \item and if $q \in Q(J_k)$ and $B \in \Below(J_k)$, then
        \[
            q \odot B = \bullet  
                \quad \text{and} \quad
            \prodout{q}{B} = \prodval{q} \cdot B.
        \]
\end{itemize}


\subsection{The regular case}
Suppose that $J_k \in \Core(\cK ; \beta)$ is regular.

\subsubsection{States and values}
The state set for a regular $\GJ$-class $J_k$ is 
    \[
        Q(J_k) = J_k \cup \{ \bullet \}.
    \]
The value map on non-initial states is given by $\prodval{X} = X$.

\subsubsection{Action and output}
The action and output of $A \in J_k$ at $\bullet$ are given by
    \[
        \bullet \odot A = A
        \quad \text{and} \quad
        \prodout{\bullet}{A} = I
    \]
and at $X \in J_k$ by
    \[
        X \odot A = 
            \begin{cases}
            XA, 
                & \text{if $XA \gR X$;} \\
            A, 
                & \text{if $XA \slgR X$;}
            \end{cases}
                                    \quad \text{and} \quad
        \prodout{X}{A} =
            \begin{cases}
            I 
                & \text{if $XA \gR X$;} \\
            XE_A 
                & \text{if $XA \slgR X$.}
            \end{cases}
    \]


\subsection{The null case}
Suppose that $J_k \in \Core(\cK ; \beta)$ is null.

\subsubsection{}
A string $(X_1, X_2, \dots, X_n) \in J_k^+$ is said to be \textbf{fragile} if the  $\lgL$-chain
    \[
        X_1 X_2 \cdots X_{n-1} X_n \;\lgL\; X_2 \cdots X_{n-1} X_n
            \;\lgL \cdots \lgL\;
        X_{n-1} X_n \;\lgL\; X_n
    \]
is strict throughout.
Clearly the length of a fragile string must be less than or equal to $d_k$, and thus the set of fragile strings over $J_k$ is finite.

A string which is not fragile is said to be \textbf{sturdy}.
Let $\mathsf{Fragile}(J_k)$ denote the set of fragile strings over $J_k$, and let 
    \[
        \mathsf{Sturdy}(J_k) =
        \left( \bigcup_{\ell = 1}^{d_k + 1} J_k^{(\ell)} \right)
        \setminus \mathsf{Fragile}(J_k).
    \]
In words, $\mathsf{Sturdy}(J_k)$ is the set of sturdy strings of length less than or equal to $d_k + 1$.
Note that all strings of length greater than or equal to $d_k + 1$ are sturdy.

\subsubsection{}
The action and output will utilize two functions
    \[
        \mathsf{cut} : \mathsf{Sturdy}(J_k) \arw \mathsf{Fragile}(J_k)
        \quad\quad \text{and} \quad\quad
        \mathsf{send} : \mathsf{Sturdy}(J_k) \arw \cB(J_k)
    \]
which act on a sturdy string $\bdx = (X_1 , \dots, X_n)$ by
    \begin{align*}
        (X_1 , \dots, X_n)\mathsf{cut} \; \quad &= \quad (X_{M+1},  \dots, X_n) \\
           &\text{and}  \\
        (X_1 , \dots, X_n)\mathsf{send} \quad &= \quad X_1 \cdots X_M E_{(X_{M+1} \cdots X_n)},
    \end{align*}
where the ``cutting point'' $M$ is defined to be
    \[
        M = \max \setst{1 \leq m \leq n-1}{X_m X_{m+1} \cdots X_n \gL X_{m+1} \cdots X_n },
    \]
which is well-defined since $\bdx$ is sturdy.

\subsubsection{States and values}
The state set for a null $\GJ$-class $J_k$ is 
    \[
        Q(J_k) = \mathsf{Fragile}(J_k) \cup \{ \bullet \}.
    \]
The value map on non-initial states is given by 
    \[
        \prodval{(X_1, \dots, X_n)} = X_1 \cdots X_n.
    \]

\subsubsection{Action and output}
The action and output of $A \in J_k$ at $\bullet$ are given by
    \[
        \bullet \odot A = (A)
        \quad \text{and} \quad
        \prodout{\bullet}{A} = I;
    \]
and at $(X_1, \dots, X_n) \in \mathsf{Fragile}(J_k)$ by
    \[
        (X_1, \dots, X_n) \odot A = 
            \begin{cases}
            (X_1, \dots, X_n, A), 
                & \text{if $(X_1, \dots, X_n, A)$ is fragile;} \\
            (X_1, \dots, X_n, A)\mathsf{cut},
                & \text{otherwise;}
            \end{cases}
    \]
and
    \[
        \prodout{(X_1, \dots, X_n)}{A} =
            \begin{cases}
            I 
                & \text{if $(X_1, \dots, X_n, A)$ is fragile;} \\
            (X_1, \dots, X_n, A)\mathsf{send} 
                & \text{otherwise.}
            \end{cases}
    \]


\subsection{Global construction}    \label{prodexp:globalaction}
We will now combine the local components defined above.
For each $s \in S$ and each $\bdq = [q_N, \dots, q_1] \in Q(J_N) \times \cdots \times Q(J_1)$, define
    \[
        (\bdq)\partial_1^s = \clos{\{s\}}
        \;\; \text{and} \;\;
        (\bdq)\partial_{k+1}^s = 
            \begin{cases}
            \clos{\prodout{q_k}{(\bdq)\partial_k^s}},
                & \text{if $\prodout{q_k}{(\bdq)\partial_k^s} \neq I$; } \\
            I, 
                & \text{otherwise.}
            \end{cases}
    \]
With this, define
    \[
        \bdq \pxact s
            = 
        \bigg[\big(q_N \odot (\bdq)\partial_N^s\big), \; 
        \big(q_{N-1} \odot (\bdq)\partial_{N-1}^s\big), \; 
            \dots\; , \;
        \big(q_1 \odot (\bdq)\partial_1^s\big)\bigg].
    \]

\begin{notation}
\hfill
\begin{enumerate}
    \item Extend the notation of the action defined in \ref{prodexp:globalaction} to strings $\sigma \in S^+$ in the usual manner; e.g., if $\sigma = (s_1 s_2 \dots s_n)$, then
        \[
            \bdq \pxact \sigma \quad = \quad ( \cdots ((\bdq \pxact s_1) \pxact s_2 ) \cdots ) \pxact s_n.
        \]
    \item Denote the state for whom all coordinates are $\bullet$ by $[\bullet]$; i.e.,
        \[
            [\bullet] = [\bullet, \bullet, \dots, \bullet].
        \]
\end{enumerate}
\end{notation}

\begin{defn}
The \textbf{product expansion} of $S$ with respect to the setup data $\mathfrak{d}$ is the automaton
    \[
        \PRXP(S ; \mathfrak{d}) 
            = 
        (S, \;\STATES{S}{\mathfrak{d}},\; [\bullet], \pxact)
    \]
where $[\bullet]$ and $\pxact$ are as above, and 
    \[
        \STATES{S}{\mathfrak{d}}
            =
        \setst{[\bullet] \pxact \sigma}{\sigma \in S^+} \cup \big\{ [\bullet] \big\}.
    \]
\end{defn}

\begin{lem} \label{prodexp:updatemult}
Let $\bdq = [q_N , q_{N-1}, \dots, q_1]$ be a state of $\PRXP(S ; \mathfrak{d})$ and let $s \in S$. Then
        \[
            \prodval{q_k} \cdot (\bdq)\partial_k^s \quad = \quad
            \prodout{q_k}{(\bdq)\partial_k^s} \cdot \prodval{ q_k \odot (\bdq)\partial_k^s }
        \]
for all coordinates $N \geq k \geq 1$.\footnotemark
\footnotetext{In words, Lemma \ref{prodexp:updatemult} says that "the (value of the) old coordinate times the input is equal to the 'raw' output times the (value of the) updated coordinate". This 'raw' output is then "intercepted and blown up" before it acts on the coordinates below it, and the axioms for blowup operators ensure that this remains well-behaved.}
\end{lem}

\begin{proof}
This follows from a routine case-checking argument.
\end{proof}

\subsection{Flow}
The product expansion $\PRXP(S ; \mathfrak{d})$ can be considered as a flow automaton when equipped with the \textbf{canonical flow}
    \[
        \flow{-}{\pi} : 
            \STATES{S}{\mathfrak{d}} \setminus \big\{ [\bullet] \big\} 
        \arw \cK
    \]
whose action on a non-initial state $\bdq = [q_N, \dots, q_1]$ is given by
    \[
        \flow{\bdq}{\pi} 
            =
        \prodval{q_N} \cdot \prodval{q_{N-1}} 
            \cdot \;\cdots\; \cdot 
        \prodval{q_2} \cdot \prodval{q_1}.
    \]

\begin{lem} \label{prodexp:canflowdefined}
The canonical flow $\flow{-}{\pi}$ is a flow whose image is contained in $\cK$.
\end{lem}

\begin{proof}
Since the image of $\flow{-}{\pi}$ is obviously contained in $\cK$, it remains only to show that it is a flow. 
Note that Lemma \ref{prodexp:updatemult} will be used a great many times throughout the proof, and so in the interest of reducing clutter we will mention its use only here.

To begin, let $s \in S$ and observe that $\flow{[\bullet] \pxact s}{\pi} = \clos{\{s\}}$, of which $s$ is clearly an element.
Now, let $\bdq = [ q_N, q_{N-1}, \dots, q_1]$ be a non-initial state and let
    \[
        A_k = 
                \begin{cases}
                    \prodout{q_{k-1}}{(\bdq) \partial_{k-1}^s}, 
                            & N+1 \geq k \geq 2; \\
                    \{s\}, 
                            & k = 1.
                \end{cases}
    \]
Note that $\clos{A_k} = (\bdq) \partial_{k}^s$ for $N \geq k \geq 1$.

Proceed by induction.
First, since $A_{N+1} = I$
    \[
        \prodval{q_N \odot \clos{A_N}}
        \;\;\; = \;\;\;
        \prodval{q_N} \cdot \clos{A_N}.
    \]
Hence
    \begin{align*}
        \prodval{q_N \odot \clos{A_N}} \cdot \prodval{q_{N-1} \odot \clos{A_{N-1}}}
        \;\; &= \;\;
        \prodval{q_N} \cdot \clos{A_N} \cdot \prodval{q_{N-1} \odot \clos{A_{N-1}}}
        \\ \\
        \;\; &\supseteq\;\; 
        \prodval{q_N} \cdot A_N \cdot \prodval{q_{N-1} \odot \clos{A_{N-1}}};
    \end{align*}
and so
    \begin{align*}
        \prodval{q_N} \cdot A_N \cdot \prodval{q_{N-1} \odot \clos{A_{N-1}}}
        \;\; &= \;\;
        \prodval{q_N} \cdot \prodval{q_{N-1}} \cdot \clos{A_{N-1}}
        \\ \\
        \;\; &\supseteq \;\;
        \prodval{q_N} \cdot \prodval{q_{N-1}} \cdot A_{N-1}.
    \end{align*}
Proceeding with the induction, if
    \[
        \left( \prod_{i=N}^k \prodval{q_i} \right) 
            \cdot \clos{A_k} 
                \;\;\; \subseteq \;\;\;
        \prod_{i=N}^k \prodval{q_i \odot \clos{A_i}},
    \]
then
\begin{align*}
    \left( \prod_{i=N}^{k-1} \prodval{q_i} \right) \cdot \clos{A_{k-1}} 
    \;\;\; &\subseteq \;\;\;
    \left( \prod_{i=N}^k \prodval{q_i} \right) 
        \cdot \clos{A_k} 
        \cdot \prodval{q_{k-1} \odot \clos{A_{k-1}}}
    \\ \\
    \;\;\; &\subseteq \;\;\;
    \prod_{i=N}^{k-1} \prodval{q_i \odot \clos{A_i}}.
\end{align*}
Thus $\flow{\bdq}{\pi} \cdot \{s\} \; \subseteq \; \flow{\bdq \pxact s}{\pi}$, as desired.
\end{proof}


%

\section{Cascades and the Zeiger property}      \label{section:cascades}

\subsection{Transformation semigroups}
A (finite) \textbf{transformation semigroup} is a pair $(X, S)$ where $X$ is a finite set and $S$ is a finite semigroup acting on the right of $X$.

\subsection{Wreath products}
Given transformation semigroups $(X , S)$ and $(Y , T)$, their \textbf{wreath product} is the transformation semigroup
    \[
        (X , S) \wr (Y , T) 
            =
        (X \times Y , S^Y \rtimes T)
    \]
where $S^Y$ denotes the semigroup of set functions from $Y$ to $S$ under pointwise multiplication, and where the action of $t \in T$ sends $f \in S^Y$ to the function given by
    \[
        \prescript{t}{}{f} : y \longmapsto (y \cdot t)f
    \]
for all $y \in Y$.
The action of this semidirect product on $X \times Y$ is given by
    \[
        (x , y) \cdot (f , t)
            =
        (x \cdot (y)f , y \cdot t)
    \]
for all $(x , y) \in X \times Y$ and all $(f , t) \in S^Y \rtimes T$.

The wreath product of transformation semigroups is associative; that is,
    \[
        \big( (X_3 , S_3) \wr (X_2 , S_2) \big) \wr (X_1, S_1)
            =
        (X_3 , S_3) \wr \big( (X_2 , S_2) \wr (X_1, S_1) \big)
    \]
for all triples of transformation semigroups $(X_i , S_i)$ ($i = 1 , 2 , 3$).

\subsection{Cascades}   \label{cascadesdef}
An element of the action semigroup of an iterated wreath product
    \[
        (X_n , S_n) \wr (X_{n-1} , S_{n-1}) \wr \cdots \wr (X_1 , S_1)
    \]
may be described as a \textbf{cascade}, which is a tuple $d = (d_n , d_{n-1} , \dots , d_1)$ with
    \begin{align*}
        d_1 & \in S_1,
        \\
        d_2 & : Q_1 \arw S_2,
        \\
        d_3 & : Q_2 \times Q_1 \arw S_3,
        \\
        &\vdots 
        \\
        d_{n-1} &: Q_{n-2} \times \cdots \times Q_2 \times Q_1 \arw S_{n-1},
        \\
        d_n &: Q_{n-1} \times Q_{n-2} \times \cdots \times Q_2 \times Q_1 \arw S_n
    \end{align*}
where $d_k$ is a set function for each $2 \leq k \leq n$.
The action of $d$ is defined in the evident coordinatewise manner.
Explicitly, given $\bdx = (x_n , x_{n-1} , \dots, x_1)$ with $x_k \in X_k$, the $i$th coordinate of $\bdx \cdot d$ is given by
    \[
        x_i \cdot (x_{i-1}, \dots , x_1)d_i
    \]
for $i \neq 1$, and the definition in the case of $i = 1$ is obvious.

\begin{notation}
Retaining notation from \ref{cascadesdef}, for convenience we write
    \[
        (\bdx)d_i = (x_{i-1} , \dots , x_1)d_i
    \]
for all $n \geq i \geq 1$.
\end{notation}

\begin{defn}
A \textbf{PR-transformation semigroup} (\underline{P}ermutation \underline{R}eset) is a transformation semigroup $(X, S)$ such that $S = G \cup C_X$ where $G$ is a finite group for which $1_G$ acts as the identity on $X$ and where $C_X$ denotes the right-zero semigroup of constant maps on $X$.
\end{defn}

\subsection{Zeiger cascades}    \label{cascades:zeigermaps}
Consider a wreath product
    \[
        (\widetilde{X} , \widetilde{S}) 
            =
        (X_n , S_n) \wr (X_{n-1} , S_{n-1}) \wr \cdots \wr (X_1 , S_1)
    \]
wherein each $(X_i , S_i)$ is a PR-transformation semigroup.
Write $S_i = G_i \cup C_{X_i}$ for each $n \geq i \geq 1$.
A cascade $d = (d_n , \dots , d_1) \in \widetilde{S}$ is said to be \textbf{Zeiger} if
    \[
        (\bdx)d_i \in G_i 
            \quad \Longrightarrow \quad
        \text{$(\bdx)d_k = 1_{G_k}$ for all $n \geq k > i$}
    \]
for every $\bdx \in \widetilde{X}$.

\begin{defn}
Retaining notation from \ref{cascades:zeigermaps}, let
    \[
        \mathcal{Z}(\widetilde{X} , \widetilde{S})
            =
        \setst{ d \in \widetilde{S} }{\text{$d$ is Zeiger}}.
    \]
\end{defn}

\begin{lem}     \label{cascades:lem:subgroups}
Any subgroup of $\mathcal{Z}(\widetilde{X} , \widetilde{S})$ is isomorphic to a subgroup of some $G_k$.
\end{lem}

\begin{proof}
The claim is obvious in the case of trivial subgroups, so we assume otherwise for hypothetical groups mentioned hereafter.
For each $n \geq k \geq 1$, let
    \[
        (\widetilde{X}_k , \widetilde{S}_k)
            =
        (X_n , S_n) \wr \cdots \wr (X_k , S_k)
    \]
and proceed by induction on $k$.

Clearly any subgroup of $(\widetilde{X}_n , \widetilde{S}_n) = (X_n, S_n)$ is a subgroup of $G_n$, so suppose that the claim holds through $k$ and let $G$ be a subgroup of $(\widetilde{X}_{k-1} , \widetilde{S}_{k-1})$.

If $e = (e_n , \dots, e_{k-1})$ is the identity of $G$, then $e_{k-1}$ is either the identity of $G_{k-1}$ or a constant map.
If $e_{k-1} = 1_{G_{k-1}}$, then $(\bdx)e_i = 1_{G_i}$ for all $\bdx \in \widetilde{X}_{k-1}$ and all $n \geq i \geq k-1$.
It follows that no element of $G$ can act as a constant map on $X_{k-1}$, and thus $(\widetilde{X}_{k-1} , G)$ is a sub-transformation group of
    \[
        (X_n , \{ 1_{G_n} \}) \wr \cdots \wr (X_k , \{ 1_{G_k} \}) \wr (X_{k-1} , G_{k-1})
    \]
whose action group is isomorphic to $G_{k-1}$.

If $e_{k-1}$ is a constant map, then it is easy to see that the $k-1$th component of any $g \in G$ must be $e_{k-1}$ as well.
Consequently the $k-1$th coordinate action of $G$ is trivial, which implies that $G$ is isomorphic to a subgroup of $(\widetilde{X}_k , \widetilde{S}_k)$ and is therefore vulnerable to the inductive hypothesis, thus yielding the lemma.\footnotemark
\footnotetext{This proof is essentially taken from \cite{GENERALHOLONOMY}.}
\end{proof}

%

%

\section{Subgroups of the product expansion}      \label{section:properties}

\subsection{}
Fix a semigroup $S$ alongside setup data $\mathfrak{d}$, and retain notation irom \ref{prodexp:notationfixed}.

\begin{notation}
\hfill
\begin{enumerate}
    \item Let $\LPXS_k$ denote the transition semigroup of $(\Inputs(J_k), Q(J_k), \bullet, \odot)$ for each $N \geq k \geq 1$; and
    \item let $\PXS$ denote the transition semigroup of the automaton $\PRXP(S ; \mathfrak{d})$.
\end{enumerate}
\end{notation}

\begin{lem}     \label{properties:lem:zeigerdiv}
There is a wreath product
    \[
        (\widetilde{P} , \widetilde{T})
            =
        (P_m , T_m) \wr (P_{m-1} , T_{m-1}) \wr \cdots \wr (P_1 , T_1)
    \]
where each $(P_i , T_i)$ is a PR-transformation semigroup with $T_i = G_i \cup C_{P_i}$ such that
\begin{enumerate}
    \item each $G_i$ is either trivial or isomorphic to a maximal subgroup of some $\beta$-core $\GJ$-class of $\cK$, and such that
    \item there exists an embedding of transformation semigroups
        \[
            \big( \STATES{S}{\mathfrak{d}} , \PXS \big) 
                \longinj 
            \big( \widetilde{P} , \mathcal{Z}(\widetilde{P} , \widetilde{T}) \big).
        \]
\end{enumerate}
\end{lem}

\begin{proof}
The basic idea for the proof is as follows.
For each $J_k \in \Core(\cK ; \beta)$, we define a wreath product
    \[
        (\widetilde{P}_k , \widetilde{T}_k)
            =
        (P_{k , m_k} , T_{k , m_k}) \wr \cdots \wr (P_{k , 1} , T_{k , 1})
    \]
along with an embedding of transformation semigroups
    \[
        (\chi_k , \delta_k) :
            (Q(J_k) , \LPXS_k)
        \longinj
            (\widetilde{P}_k , \mathcal{Z}(\widetilde{P}_k , \widetilde{T}_k)).
    \]
For each $J_k$, let $\Omega_k = \Img(\chi_k)$ and let $Z_k = \Img(\delta_k)$.
We then consider 
    \[
        (\widetilde{\Omega} , \widetilde{Z})
            =
        (\Omega_N, Z_N) \wr (\Omega_{N-1} , Z_{N-1}) \wr \cdots \wr (\Omega_1 , Z_1)
    \]
where, as before, $N$ denotes the cardinality of $\Core(\cK ; \beta)$.
Immediately we have an embedding
    \[
        \widetilde{\chi} 
            = 
        \prod_{k=N}^1 \chi_k 
            : \STATES{S}{\mathfrak{d}} \longinj \widetilde{\Omega}.
    \]
Moreover, recalling from \ref{prodexp:globalaction} that 
    \[
        \bdq \pxact s
            = 
        \bigg[\big(q_N \odot (\bdq)\partial_N^s \big), \; 
        \big(q_{N-1} \odot (\bdq)\partial_{N-1}^s \big), \; 
            \dots\; , \;
        \big(q_1 \odot (\bdq)\partial_1^s \big)\bigg].
    \]
for all $\bdq \in \STATES{S}{\mathfrak{d}}$ and all $s \in S$, we obtain a map
    \[
        \widetilde{\delta} : \PXS \longinj \widetilde{Z}
            \quad \text{written} \quad
        (-) \pxact s    
            \longmapsto     
        \widetilde{\delta}^s = (\widetilde{\delta}^s_N , \dots , \widetilde{\delta}^s_1)
    \]
where $\widetilde{\delta}^s_1 = (\partial_1^s)\delta_1$ and
    \[
        (p_k , \cdots , p_1)\widetilde{\delta}^s_{k+1}
            =
        \Big(  \big(\bullet , \dots, \bullet, (p_k)\chi_k^\inv ,  \dots , (p_1)\chi_1^\inv\big)\partial_{k+1}^s \Big)\delta_{k+1}
    \]
for all $(p_k , \cdots , p_1) \in \prod_{i=k}^1 \Omega_i$ once $\widetilde{\delta}^s_k$ is defined.

Since each $\widetilde{\delta}^s_k$ is Zeiger on $(\Omega_k , Z_k)$, the observation that any ``local'' group action yields $I$ as output implies that each $\widetilde{\delta}^s$ is ``globally'' Zeiger as well.
Once established, this yields the lemma.

Within the various constructions it will be convenient to write $X^\bullet = X \cup \{ \bullet \}$ (where $X$ is a set and the union is disjoint).
As before, the details of the construction for each $J_k$ depend on whether $J_k$ is regular or null.

First, suppose that $J_k \in \Core(\cK ; \beta)$ is regular.
Fix a coordinatization
    \[
        \lambda : J_k \arw A \times G \times B
            \quad \text{written} \quad
        X \longmapsto (a_X, g_X, b_X),
    \]
where $A$ and $B$ respectively index the $\GR$- and $\GL$-classes of $J_k$ and $G$ is the Sch{\"u}tzenberger group of $J_k$, along with a structure matrix
    \[
        C = \strmat{- , -} : B \times A \arw G^0
    \]
such that 
    \[
        (X)\lambda \cdot (Y)\lambda 
            =
        (a_X , \, g_X \cdot \strmat{b_X , a_Y} \cdot g_Y , \, b_Y)
    \]
for any $X , Y \in J_k$ with $\strmat{b_X , a_Y} \neq 0$.
The wreath product in this case is
    \[
        (\widetilde{P}_k , \widetilde{T}_k)
            =
        \left( (A \times G)^\bullet , \; G \cup C_{(A \times G)^\bullet} \right)
            \wr
        \Big( B^\bullet , \; \{ \id{{B^\bullet}} \} \cup C_{B^\bullet} \Big),
    \]
where the action of $G$ on $(A \times G)^\bullet$ is given by
    \[
        \bullet \cdot g = \bullet
            \quad \text{and} \quad
        (a , h) \cdot g = (a , hg)
    \]
The embedding $\chi_k : Q(J_k) \inj \widetilde{P}_k$ is given by
    \[
        (\bullet)\chi_k = (\bullet, \bullet)
            \quad \text{and} \quad
        (X)\chi_k = ( ( a_X , g_X ) , b_X ).
    \]
The embedding $\delta_k : \LPXS_k \inj \widetilde{T}_k$ sends the identity and ``reset-to-$\bullet$'' maps to the evident elements of $\widetilde{T}_k$; and if $X \in J_k$ then $(-) \odot X$ is sent to the cascade
    \[
        \delta^X_k 
            =
        \left( \delta^X_{k , 2} , \delta^X_{k , 1} \right)
    \]
where $\delta^X_{k , 1}$ is the constant map $b \mapsto b_X$, and the actions assigned by $\delta^X_{k , 2}$ are defined at $\bullet$ by 
    \[
        \bullet \cdot (b)\delta^X_{k , 2}
            \quad=\quad
        \bullet \cdot (\bullet)\delta^X_{k , 2}
            \quad=\quad
        (a_X , g_X);
    \]
and at $(a , g) \in A \times G$ by
    \[
        ( a , g ) \cdot (\bullet)\delta^X_{k , 2}
            =
        ( a_X , g_X )
    \]
when given $\bullet$, and by
    \[
        ( a , g ) \cdot (b)\delta^X_{k , 2}
            =
        \begin{cases}
            ( a  ,  g \cdot \strmat{b , a_X} \cdot g_X ),    
                    &
                    \text{if $\strmat{b , a_X} \neq 0$};
        \\
            ( a_X , g_X ),
                    &
                    \text{if $\strmat{b , a_X} = 0$}.
        \end{cases}
    \]
when given $b \in B$.
Each such cascade is easily seen to be Zeiger, thus giving us permission to move onto the null case.

When $J_k$ is null, the wreath product is given by
    \[
        (\widetilde{P}_k , \widetilde{T}_k)
            =
    \underbrace{
        \left( J_k^\bullet  \, , \; \Set*{ \id{{J_k^\bullet}} } \cup C_{{J_k^\bullet}} \right)
            \wr
        \left( J_k^\bullet  \, , \;  \Set*{ \id{{J_k^\bullet}} } \cup C_{{J_k^\bullet}} \right)
            \wr \cdots \wr 
        \left( J_k^\bullet  \, , \;  \Set*{ \id{{J_k^\bullet}} } \cup C_{{J_k^\bullet}} \right)
    }_{\text{$d_k + 1$ times}}
    \]
where $d_k$ is defined as in \ref{prodexp:notationfixed}.
The embedding $\chi_k : Q(J_k) \inj \widetilde{P}_k$ is given by
    \[
        (\bullet)\chi_k = (\bullet, \dots , \bullet)
            \quad \text{and} \quad
        (X_1 , \dots , X_n)\chi_k = (\bullet , \dots , \bullet , X_1 , \dots , X_n).
    \]
Defining the embedding $\delta_k : \LPXS_k \inj \widetilde{T}_k$---which is written as
    \[
        (-) \odot \sigma 
            \quad \longmapsto \quad
        \delta^\sigma_k 
            =
        \left( \delta^\sigma_{k , {d_k}+1} , \delta^\sigma_{k , {d_k}} , \dots , \delta^\sigma_{k , 1} \right)
    \]
for any $\sigma \in \Inputs(J_k)$---is slightly more involved.
The cascades assigned to the actions of $I$ and $X \in \Below(J_k)$ are respectively defined by
    \[
        \bdx \cdot \delta^I_k  =  \bdx
            \quad \text{and} \quad
        \bdx \cdot \delta^X_k  =  (\bullet, \dots , \bullet)
    \]
at every $\bdx \in \widetilde{P}_k$.
To define $\delta^X_k$ in remaining case where $X \in J_k$, we must first define a ``condensation'' map
    \[
        \mu : \widetilde{P}_k \arw \left(  \bigcup_{m=1}^{d_k + 1} J_k^{(m)} \right)^\bullet
    \]
which sends $(\bullet , \dots , \bullet)$  to $(\bullet)$ and which otherwise ``deletes''  occurrences of $\bullet$ appearing in a string alongside members of $J_k$.
Next, define a function
    \[
        M_\ast : \widetilde{P}_k \arw \{ 1 , 2 , \dots , d_k + 1 \} \cup \{ \infty \}
    \]
by sending $( X_{{d_k} + 1} , \dots , X_1 ) \in \widetilde{P}_k$ to the value 
    \[
        M_{( X_{{d_k} + 1} , \dots , X_1 )}
            =
        \begin{cases}
            \min \setst{ m \in \{ 1 , \dots , d_k + 1 \} }{\text{$(X_m , \dots , X_1)\mu$ is sturdy}},
                &
                \text{if it exists};
        \\
            \infty,
                &
                \text{otherwise}.
        \end{cases}
    \]
Now, given $X \in J_k$, write the action of $\delta^X_k$ on $( X_{{d_k} + 1} , \dots , X_1 ) \in \widetilde{P}_k$ as
    \[
        ( X_{{d_k} + 1} , \dots , X_1 ) \cdot \delta^X_k
            =
        ( Y_{{d_k} + 1} , \dots , Y_1 ) 
    \]
so that we may set
    \[
        Y_i =
            \begin{cases}
                X,          &   i = 1;   
            \\
                X_{i-1},    &   2 \leq i \leq \min\Set*{ M_{( X_{{d_k} + 1} , \dots , X_1 )} , d_k + 1 };
            \\
                \bullet,    &   \text{else}
            \end{cases}
    \]
for each coordinate $Y_i$ with $i \in \{ 1 , \dots , d_k + 1 \}$.
It is once again straightforward to show that all cascades defined in this case are Zeiger.

Having covered both the null and regular cases, the promises made at the outset of this proof have been fulfilled.
From here, the tedious but straightforward process described above may be carried out to finally yield the lemma.
\end{proof}

\begin{prop}    \label{properties:prop:subgroups}
Every subgroup of $\PXS$ is isomorphic to some subgroup of a $\beta$-core $\GJ$-class of $\cK$.
\end{prop}

\begin{proof}
This follows from Lemmas \ref{properties:lem:zeigerdiv} and \ref{cascades:lem:subgroups}.
\end{proof}

%

%

\section{Upper bound and main result}      \label{section:upperbound}

\subsection{}
Recall from Section \ref{section:lowerbound} that $\CGH(S) \subseteq \PGH(S)$ for every finite semigroup $S$, where $\CGH(S)$ is the minimal $S$-complex satisfying the condition that
    \[
        \mathscr{X} \in \PL{\pvH}(\mathcal{G})
            \quad \Longrightarrow \quad
        \bigcup \mathscr{X} \in \CGH(S)
    \]
whenever $\mathcal{G}$ is a subgroup of $\CGH(S)$.
In short, $\CGH$ provides a lower bound for $\pvGH$-pointlikes.
In this section, we use the product expansion to show that this bound is exact.

\begin{notation}
For each $X \in \CGH(S)$, let $G_X$ denote the right Sch{\"u}tzenberger group of the $\GL$-class of $X$.
\end{notation}

\subsection{Blowup operator}
Define a function $\blgh : \CGH(S) \rightarrow \CGH(S)$ by setting
    \[
        (X)\blgh 
            = 
        \bigcup \setst{ X \cdot g }{ g \in K_{\pvH , {G_X}}}
    \]
for all $X \in \CGH(S)$.

\begin{lem} \label{ub:lem:blowupcore}
The map $\blgh$ is a preblowup operator on $\CGH(S)$, and the $\blgh$-core $\GJ$-classes of $\CGH(S)$ are precisely those whose Sch{\"u}tzenberger groups belong to $\pvH$.
\end{lem}

\begin{proof}
The second claim is obvious so long as the first claim holds, and so it is on the first claim that we set our sights.
Three conditions must be verified.
\begin{enumerate}
    \item Given $X \in \CGH(S)$, the fact that $1_{G_X} \in K_{\pvH , {G_X}}$ implies that $X \subseteq (X)\blgh$.
    \item If $X , Y \in \CGH(S)$ and $X \gL Y$, then $G_X = G_Y$.
    Since $\GL$ is stable under right multiplication, it follows that $(X)\blgh \gL (Y)\blgh$.
    \item Finally, given $X \in \CGH(S)$, we must show that either $(X)\blgh = X$ or $(X)\blgh \slgH X$---i.e., under the action of $\blgh$, either $X$ is fixed or $X$ strictly falls in the $\GH$-order.
    
    To begin, recall that there exists a subgroup $\widetilde{G}$ of the right stabilizer of $\LCL{X}$ which maps onto $G_X$ under the natural map.
    Hence
        \begin{align*}
            (X)\blgh 
                &=
            \bigcup \setst{ X \cdot g }{ g \in K_{\pvH , {G_X}}}
        \\      &= 
            \bigcup \setst{ X \cdot \tilde{g} }{ \tilde{g} \in K_{\pvH , \widetilde{G}}}    
        \\      &=
            X \cdot \bigcup K_{\pvH , \widetilde{G}},
        \end{align*}
    and thus $(X)\blgh \lgR X$.
    The dual argument involving \textit{left} Sch{\"u}tzenberger groups yields the dual inequality---concretely, that $(X)\blgh \lgL X$---from which it follows that $(X)\blgh \lgH X$.
    
    No effort is required in the case where $(X)\blgh \slgH X$, so suppose that $(X)\blgh \gH X$. 
    If $k \in K_{\pvH , {G_X}}$ then
        \[
            (X)\blgh \cdot k
                =
            \bigcup \setst{ X \cdot gk }{ g \in K_{\pvH , {G_X}} }
                =
            (X)\blgh,
        \]
    which implies (by faithfulness of the $G_X$-action) that $K_{\pvH , {G_X}}$ is trivial, which in turn implies that $(X)\blgh = X$.
\end{enumerate}
\end{proof}

\subsection{Setup data}
Define a tuple of setup data 
    \[
        \mathfrak{d}_{\pvGH}
            =
        \left( \CGH(S) , \beta_\pvGH^\omega , E_\ast , J_\ast \right),
    \]
where the specificities of $E_\ast$ and $J_\ast$ are arbitrary.

\begin{notation}
Denote the product expansion with respect to $\mathfrak{d}_{\pvGH}$ by
    \[
        \PRXP(S ; \pvGH) 
            =
        (S, \;\STATES{S}{\pvGH},\; [\bullet], \pxact),
    \]
and denote its associated canonical flow by
    \[
        \flow{-}{\pvGH} : 
            \STATES{S}{\pvGH} \setminus \big\{ [\bullet] \big\} 
        \arw \CGH(S).
    \]
Moreover, let $\PXSGH$ denote the transition semigroup of $\PRXP(S ; \pvGH)$.
\end{notation}

\begin{prop}    \label{upperbound:prop:ub}
Retaining notation from above,
    \[
        \Cov \left( S , \PRXP(S ; \pvGH) , \flow{-}{\pvGH} \right) 
            \subseteq 
        \CGH(S)
    \]
and $\PXSGH$ is a member of $\pvGH$.
Consequently, $\PGH(S) \subseteq \CGH(S)$.
\end{prop}

\begin{proof}
The first claim is a particular case of Lemma \ref{prodexp:canflowdefined}, and the second claim follows from Lemma \ref{ub:lem:blowupcore} by way of Proposition \ref{properties:prop:subgroups}.
The conclusion then follows from Proposition \ref{flowequivalence}.
\end{proof}

\begin{thm}
If $\pvH$ is a pseudovariety of groups, then $\CGH = \PGH$.
\end{thm}

\begin{proof}
Combine \ref{lowerbound:lb} and Proposition \ref{upperbound:prop:ub}.
\end{proof}

\begin{cor}
If $\pvH$ is a decidable pseudovariety of groups, then $\PGH$ is decidable.
\end{cor}

%

\printbibliography
\hrulefill
\end{document}